\documentclass[12pt]{article}

\usepackage{amssymb}
\usepackage{amsfonts}
\usepackage{amsmath}
\usepackage{amsthm}
\usepackage[margin=0.80in]{geometry}
\usepackage{enumerate}
\usepackage{amstext}
\usepackage{layout}

\numberwithin{equation}{section}

\newtheorem{theorem}{Theorem}[section]

\newtheorem{lemma}{Lemma}[section]

\newtheorem{remark}{Remark}[section]

\renewcommand{\thefootnote}{\arabic{footnote}}

\begin{document}

\noindent {\bf\large{A supplement to the laws of large numbers and the large deviations}}

\vskip 0.3cm

	\noindent {\bf Deli Li\footnote{Deli Li, Departments of Mathematical Sciences, Lakehead University, Thunder Bay, Ontario P7B 5E1, Canada~~
			e-mail: dli@lakeheadu.ca} $\cdot$ Yu Miao\footnote{Yu Miao, College of Mathematics and Information Science, 
			Henan Normal University, Xinxiang, Henan 453007, China~~
			e-mail: yumiao728@gmail.com} $^*$\let\thefootnote\relax\footnote{$^*$Corresponding author: Yu Miao (Telephone:
			+86-373-332-8615, FAX: +86-373-332-8615)}}

\vskip 0.3cm
	
\noindent {\bf Abstract}~~Let $0 < p < 2$. Let $\{X, X_{n}; n \geq 1\}$ be a sequence of independent and identically 
distributed $\mathbf{B}$-valued random variables and set $S_{n} = \sum_{i=1}^{n}X_{i},~n \geq 1$. 
In this paper, a supplement to the classical laws of large numbers and the classical large deviations 
is provided. We show that if $S_{n}/n^{1/p} \rightarrow_{\mathbb{P}} 0$, then, 
for all $s > 0$,
\[
\limsup_{n \to \infty} \frac{1}{\log n} \log \mathbb{P}\left(\left\|S_{n} \right\| > s n^{1/p} \right)
= - (\bar{\beta} - p)/p
\]
and 
\[
\liminf_{n \to \infty} \frac{1}{\log n} \log \mathbb{P}\left(\left\|S_{n} \right\| > s n^{1/p} \right)
= -(\underline{\beta} - p)/p,
\]
where
\[
\bar{\beta} = - \limsup_{t \rightarrow \infty} \frac{\log \mathbb{P}(\log \|X\| > t)}{t}
~~\mbox{and}~~\underline{\beta} = - \liminf_{t \rightarrow \infty} \frac{\log \mathbb{P}(\log \|X\| > t)}{t}.
\]
The main tools employed in proving this result are the symmetrization technique and three powerful inequalities 
established by Hoffmann-J{\o}rgensen (1974), de Acosta (1981), and Ledoux and Talagrand (1991), respectively.
As a special case of this result, the main results of Hu and Nyrhinen (2004) are not only improved, but also 
extended.

~\\

\noindent {\bf Keywords}~~ Laws of large numbers $\cdot$ Large deviations $\cdot$ Heavy-tailed random variables $\cdot$ Logarithmic 
asymptotic behaviors $\cdot$ Sums of i.i.d. random variables

\vskip 0.3cm

\noindent {\bf Mathematics Subject Classification (2010)} Primary~60F10 $\cdot$ Secondary 60B12 $\cdot$ 60F05 $\cdot$ 60G50

\vskip 0.3cm

\noindent {\bf Running Head}:~~The laws of large numbers and the large deviations

\section{Introduction}

Throughout this paper, let $(\mathbf{B}, \| \cdot \| )$ be a real separable Banach space equipped with 
its Borel $\sigma$-algebra $\mathcal{B}$ ($=$ the $\sigma$-algebra generated by the class of open subsets 
of $\mathbf{B}$ determined by $\|\cdot\|$). Let $\{X, X_{n};~n \geq 1 \}$ be a sequence of independent 
and identically distributed (i.i.d.) $\mathbf{B}$-valued random variables defined on a probability space 
$(\Omega, \mathcal{F}, \mathbb{P})$. As usual, let $S_{n} = \sum_{k=1}^{n} X_{k},~ n \geq 1$ denote 
their partial sums. 

Many questions in probability theory can be formulated as a law of large numbers. To bring into focus 
of essence of this paper, we begin the statements of the classical laws of large numbers. If $0 < p < 2$ 
and if $\{X, X_{n};~n \geq 1 \}$ is a sequence of i.i.d. real-valued random variables (that is, if 
$\mathbf{B} = \mathbb{R}$), then
\begin{equation}
\lim_{n \rightarrow \infty} \frac{S_{n}}{n^{1/p}} = 0 ~~\mbox{almost surely (a.s.)~~ if and only if}~
\mathbb{E}|X|^{p} < \infty ~\mbox{where}~ \mathbb{E}X = 0
~\mbox{whenever}~ p \geq 1
\end{equation}
and
\begin{equation}
\frac{S_{n}}{n^{1/p}} \rightarrow_{\mathbb{P}} 0~~\mbox{if and only if}~~
	\left \{
	\begin{array}{ll}
	\mbox{$\displaystyle \lim_{n \rightarrow \infty} 
		n \mathbb{P}\left(|X| > n^{1/p} \right) = 0$} & \mbox{$\displaystyle \mbox{if}~ 0 < p < 1$,}\\
	&\\
	\mbox{$\displaystyle \lim_{n \rightarrow \infty} \left(\mathbb{E}XI(|X| \leq n) +
		n \mathbb{P}\left(|X| > n \right) \right) = 0$} & \mbox{$\displaystyle \mbox{if}~ p = 1$,}\\
	&\\
	\mbox{$\displaystyle \mathbb{E}X = 0 ~\mbox{and}~ \lim_{n \rightarrow \infty} 
		n \mathbb{P}\left(|X| > n^{1/p} \right) = 0$} & \mbox{$\displaystyle \mbox{if}~ 1 < p < 2$,}\\
	\end{array}
	\right.
	\end{equation}
where ``$\rightarrow_{\mathbb{P}}$" stands for convergence in probability. Statement (1.1) is the famous
Kolmogorov-Marcinkiewicz-Zygmund strong law of large numbers (SLLN) (see Kolmogoroff (1930) for $p = 1$ and
Marcinkiewicz and Zygmund (1937) for $p \neq 1$) and statement (1.2) is the celebrated Kolmogorov-Feller 
weak law of large numbers (WLLN) (see the books by Gnedenko and Kolmogorov (1968) and Feller (1971)). 

The classical Kolmogorov SLLN in real separable Banach spaces was established by Mourier (1953). The extension of the
Kolmogorov-Marcinkiewicz-Zygmund SLLN to $\mathbf{B}$-valued random variables is independently due to Azlarov and Volodin
(1981, Theorem) and de Acosta (1981, Theorem 3.1) who showed that, if $0 < p < 2$ and if $\{X, X_{n};~n \geq 1 \}$ is a sequence of i.i.d. 
$\mathbf{B}$-valued random variables, then
\[
\lim_{n \rightarrow \infty} \frac{S_{n}}{n^{1/p}} = 0~~\mbox{a.s. if and only if}~~
\mathbb{E}\|X\|^{p} < \infty~~\mbox{and}~~\frac{S_{n}}{n^{1/p}}
\rightarrow_{\mathbb{P}} 0.
\] 
In addition, de Acosta (1981, Theorem 4.1) provided a remarkable characterization of Rademacher 
type $p$ ($1 \leq p < 2$) Banach spaces via the SLLN by showing that the following two statements 
are equivalent:
\begin{align*}
& {\bf (i)} \quad \mbox{The Banach space $\mathbf{B}$ is of
	Rademacher type $p$.}\\
& {\bf (ii)} \quad \mbox{For every sequence $\{X, X_{n}; ~n \geq 1 \}$
	of i.i.d. $\mathbf{B}$-valued random variables},
\end{align*}
\[
\lim_{n \rightarrow \infty} \frac{S_{n}}{n^{1/p}} = 0~~\mbox{a.s. if
	and only if}~~\mathbb{E}\|X\|^{p} < \infty~~\mbox{and}~~\mathbb{E}X
= 0.
\] 
We refer to Ledoux and Talagrand (1991) for the definitions of Rademacher type $p$
and stable type $p$ Banach spaces. Marcus and Woyczy\'{n}ski (1979, Theorem 5.1) provided 
a remarkable characterization of stable type $p$ ($1 \leq p < 2$) 
Banach spaces via the WLLN by showing that the following two statements are equivalent:
\begin{align*}
& {\bf (i)} \quad \mbox{The Banach space $\mathbf{B}$ is of
	stable type $p$.}\\
& {\bf (ii)} \quad \mbox{For every sequence $\{X, X_{n}; ~n \geq 1 \}$
	of i.i.d. symmetric $\mathbf{B}$-valued random variables},
\end{align*}
\[
 \frac{S_{n}}{n^{1/p}} \rightarrow_{\mathbb{P}} 0~~\mbox{if
	and only if}~~\lim_{n \rightarrow \infty} n \mathbb{P}\left(\|X\| > n^{1/p} \right) = 0.
\] 

Obviously, the SLLN and WLLN results say little or nothing about the rate of convergence, however, 
which is often important for many applications of probability theory.

Let $\{X, X_{n}; n \geq 1 \}$ be a sequence of i.i.d. real-valued random variables.
Cram\'{e}r (1938) and Chernoff (1952) showed that if 
\begin{equation}
M(t) \equiv \mathbb{E} \left(e^{tX} \right) < \infty ~~\forall ~t \in \mathbb{R},
\end{equation}
then
\begin{description}
\item {\bf (i)}~ for every closed set ${\bf A} \subseteq \mathbb{R}$,
\[
\limsup_{n \to \infty} \frac{1}{n} \log \mathbb{P}\left(S_{n}/n \in {\bf A} \right) \leq - \Lambda({\bf A}),
\]
\item {\bf (ii)}~ for every open set ${\bf A} \subseteq \mathbb{R}$,
\[
\liminf_{n \to \infty} \frac{1}{n}\log \mathbb{P}\left(S_{n}/n \in {\bf A} \right) \geq - \Lambda({\bf A}),
\]
\end{description}
where, for $x \in \mathbb{R}$ and ${\bf A} \subseteq \mathbb{R}$, 
$\lambda(x) = \sup_{t \in \mathbb{R}}\left(tx - \log M(t) \right)$, 
$\Lambda({\bf A}) = \inf_{x \in {\bf A}} \lambda(x)$. This fundamental result is what we call the
large deviation principle (LDP) for partial sums $\{S_{n}, n \geq 1 \}$.  Clearly, under condition (1.3), 
the LDP characterizes the exponential concentration behavior, as $n \rightarrow \infty$, of a
sequence of probabilities $\displaystyle \left\{\mathbb{P}\left(S_{n}/n \in \mathbf{A} \right);~n \geq 1 \right\}$ 
in terms of a rate function $\Lambda(\mathbf{A})$, ${\bf A} \subseteq \mathbb{R}$. Donsker and Varadhan (1976) and Bahadur 
and Zabell (1979) established a LDP for sums of i.i.d. $\mathbf{B}$-valued random variables. The large deviations theory 
has applications in many different scientific fields, ranging from queuing theory to statistics and from finance to engineering. 

Note that (1.3) is the moment generating function of $X$ which is a very strong moment condition of $X$. Thus, one may ask, 
naturally, whether there is some LDP analog under $S_{n}/n^{1/p} \rightarrow_{\mathbb{P}} 0$ (i.e., WLLN condtion) only. 
This paper will provide a positive answer; see Theorem 2.1 in Section 2. 

The second motivation of this paper comes from the following logarithmic asymptotic behaviors obtained by Gantert (2000) 
for the partial sums of i.i.d. non-negative random variables.

\vskip 0.2cm

\begin{theorem}
	{\rm (Gantert (2000, Theorem 2))} Let $\{X, X_{n}; n \geq 1\}$ be a sequence of non-negative i.i.d. random variables 
	and set $S_{n} = \sum_{i=1}^{n}X_{i},~n \geq 1$. Let $\alpha > 1$. Then the following holds.
	\begin{description}
		\item {\rm (a)} If $\mathbb{E}X^{\alpha} = \infty$ and $\mathbb{P}(X > t) \geq L(t)/t^{\alpha}$ for some slowly 
		varying function $L$, then for every $m > 0$ and $x \geq 1$,
		\[
		\liminf_{n \rightarrow \infty} \frac{1}{\log n} \log \mathbb{P}\left(S_{n} \geq n^{x} m \right) \geq - \alpha x + 1.
		\]
		\item {\rm (b)} If $\mathbb{E}X^{p} < \infty$ for each $p < \alpha$, then for
		every $m > \mathbb{E}X$ and $x \geq 1$,
		\[
		\limsup_{n \rightarrow \infty} \frac{1}{\log n} \log \mathbb{P}\left(S_{n} \geq n^{x} m \right) \leq - \alpha x + 1. 
		\]
    \end{description}
\end{theorem}

To further investigate the logarithmic asymptotic behaviors for the partial sums $S_{n}$, Hu and Nyrhinen (2004) introduced the 
following two parameters for non-negative random variable $X$, namely
\[
\bar{\alpha} = - \limsup_{t \rightarrow \infty} \frac{1}{t} \log \mathbb{P}(\log X > t) \in [0, \infty]
\]
and
\[
\underline{\alpha} = - \liminf_{t \rightarrow \infty} \frac{1}{t} \log \mathbb{P}(\log X > t) \in [0, \infty]. 
\]
Clearly, $\bar{\alpha} \leq \underline{\alpha}$. As pointed out by Hu and Nyrhinen (2004), a useful fact is that
\[
\bar{\alpha} = \sup\left\{\lambda \geq 0;~ \mathbb{E}X^{\lambda} < \infty \right\},
\]
and hence, if $\bar{\alpha} < \infty$, then $X$ is heavy tailed, namely, $\mathbb{E}e^{\lambda X} = \infty$ for
every $\lambda >0$. 

Write $\bar{x} = \max\{1, 1/\bar{\alpha} \}$ if $\bar{\alpha} > 0$ and $\underline{x} = \max\{1, 1/\underline{\alpha} \}$.
Hu and Nyrhinen (2004) established the following large deviations for the partial sums of non-negative i.i.d. random
variables.

\vskip 0.2cm

\begin{theorem}
	{\rm (Hu and Nyrhinen (2004, Theorems 2.1 and 2.2))} Let $\{X, X_{n}; n \geq 1\}$ be a sequence of non-negative i.i.d. random variables. We have:\\
	   (i)~~Assume that $\bar{\alpha} \in (0, \infty)$. Then for every $x > \bar{x}$,
		\[
		\limsup_{n \rightarrow \infty} \frac{1}{\log n} \log \mathbb{P}\left(S_{n} > n^{x} \right) = - \bar{\alpha} x + 1.
		\]
		If in addition $\underline{\alpha} < \infty$, then for every $x > \bar{x}$, 
		\[
		\liminf_{n \rightarrow \infty} \frac{1}{\log n} \log \mathbb{P}\left(S_{n} > n^{x} \right) = - \underline{\alpha} x + 1. 
		\]
		(ii)~~If $\bar{\alpha} = 0$, then for every $x > 1$,
		\[
		\limsup_{n \rightarrow \infty} \mathbb{P}\left(S_{n} > n^{x} \right) = 1.
		\]
		(iii)~~If $\bar{\alpha} = \infty$, then for every $x > 1$,
		\[
		\lim_{n \rightarrow \infty} \frac{1}{\log n} \log \mathbb{P}\left(S_{n} > n^{x} \right) = - \infty.
		\]
\end{theorem}

\vskip 0.2cm

Motivated the above-mentioned results, Miao, Xue, Wang, and Zhao (2012) established the logarithmic 
asymptotic behaviors for the cases of m-stationary sequences and stationary negatively associated sequences. 
Recently, Wang (2017) not only improved the logarithmic asymptotic results obtained by 
Miao, Xue, Wang, and Zhao (2012) for the partial sums of stationary sequences, but also generalized the results 
to classes of acceptable random variables and widely acceptable random variables. 

Inspired by the above discovery by Gantert (2000) and Hu and Nyrhinen (2004), in the current work a supplement 
to the classical laws of large numbers and the classical large deviations will be provided. We will establish 
the large deviations for the tail probabilities $\mathbb{P}\left(\|S_{n}\| > s n^{1/p} \right)$ for all $s > 0$ 
by giving the exact values for 
\[
\limsup_{n \rightarrow \infty} \frac{1}{\log n} \log \mathbb{P}\left(\|S_{n}\| > s n^{1/p} \right)~~\mbox{and}~~
\liminf_{n \rightarrow \infty} \frac{1}{\log n} \log \mathbb{P}\left(\|S_{n}\| > s n^{1/p} \right),
\]
where $ \{X, X_{n};~n \geq 1 \}$ is a sequence of i.i.d. {\bf B}-valued random variables and $0 < p < 2$.

The plan of the paper is as follows. The main results are stated in Section 2. Our main results, Theorems 2.1, 2.2, 
and 2.3 are stated in Section 2. For the case $1 \leq p < 2$, our results are new even in the case where the Banach 
space is the real line. Some preliminary lemmas are presented in Section 3 and the proofs of our main results are 
given in Section 4. The main tools employed in proving the main results are the remarkable Hoffmann-J{\o}rgensen (1974)
inequality, de Acosta (1981) inequality, the symmetrization technique, and an inequality of Ledoux and Talagrand (1991), etc. 

\section{Statement of the main results}

Let $X$ be a ${\bf B}$-valued random variable. Write
\begin{equation}
\bar{\beta} = - \limsup_{t \rightarrow \infty} \frac{1}{t} \log \mathbb{P}(\log \|X\| > t) \in [0, \infty]
\end{equation}
and
\begin{equation}
\underline{\beta} = - \liminf_{t \rightarrow \infty} \frac{1}{t} \log \mathbb{P}(\log \|X\| > t) \in [0, \infty]. 
\end{equation}
For any real numbers $x \in (0, \infty)$ and $y \in (-\infty, \infty)$, by convention, $(\pm\infty + y)/x =  \pm \infty$.

We now state our main results.

\begin{theorem}
Let $0 < p < 2$. Let $\{X, X_{n}; n \geq 1\}$ be a sequence of i.i.d. $\mathbf{B}$-valued random variables and 
set $S_{n} = \sum_{i=1}^{n}X_{i},~n \geq 1$. If
\begin{equation}
\frac{S_{n}}{n^{1/p}} \rightarrow_{\mathbb{P}} 0,
\end{equation}
then, for all $s > 0$,
\begin{equation}
\limsup_{n \to \infty} \frac{1}{\log n} \log \mathbb{P}\left(\left\|S_{n} \right\| > s n^{1/p} \right)
= - \left(\bar{\beta} - p \right)/p
\end{equation}
and 
\begin{equation}
\liminf_{n \to \infty} \frac{1}{\log n} \log \mathbb{P}\left(\left\|S_{n} \right\| > s n^{1/p} \right)
= -\left(\underline{\beta} - p \right)/p.
\end{equation}
Hence, under condition (2.3), for all $s > 0$,
\begin{equation}
\lim_{n \to \infty} \frac{1}{\log n} \log \mathbb{P}\left(\left\|S_{n} \right\| > s n^{1/p} \right)
= - \left(\hat{\beta} - p \right)/p ~~\mbox{if and only if}~~\bar{\beta} = \underline{\beta} = \hat{\beta}.
\end{equation}
\end{theorem}

\vskip 0.2cm

\begin{remark}
	Clearly, condition (2.3) is the weak law of large numbers for partial sums $S_{n}$, $n \geq 1$. 
	For given $0 < p < 1$, it is well known that (2.3) is equivalent to
	\[
	 \lim_{n \rightarrow \infty} 
	 n \mathbb{P}\left(\|X\| > n^{1/p} \right) = 0.
	\]
\end{remark} 

\vskip 0.2cm

\begin{remark}	
	Let $1 \leq p < 2$ and let $\{b_{n}; n \geq 1\}$ be an increasing sequences of positive 
	real numbers such that 
	\[
	\left\{b_{n}/n^{1/p}; ~n \geq 1 \right\}~\mbox{is a nondecreasing sequence}.
	\]
	Using a remarkable characterization of stable type $p$ Banach spaces
	provided by Marcus and Woyczy\'{n}ski (1979), Li, Liang, and Rosalsky (2018, Corollary 3.2) showed that, 
	if $\mathbf{B}$ is of stable type $p$, then for every sequence $\{X, X_{n}; ~n \geq 1 \}$ 
	of i.i.d. {\bf B}-valued random variables, 
	\begin{equation}
	\frac{S_{n}- \gamma_{n}}{b_{n}} 
	\rightarrow_{\mathbb{P}} 0 ~~\mbox{or}~~\limsup_{n \rightarrow \infty} 
	\mathbb{P} \left(\frac{\left\|S_{n}- \gamma_{n} \right\|}{b_{n}} 
	> \lambda \right) > 0 ~~\forall~\lambda > 0
	\end{equation}
	according as
	\begin{equation}
	\lim_{n \rightarrow \infty}
	n \mathbb{P}(\|X\| > b_{n}) = 0~~\mbox{or}~~\limsup_{n \rightarrow \infty}
	n \mathbb{P}(\|X\| > b_{n}) >0,
	\end{equation}
	where $\gamma_{n} = n \mathbb{E}\left(XI\{\|X\| \leq b_{n} \} \right)$, $n \geq 1$. 
	For the special case $b_{n} = n^{1/p}$, it can be easily deduced from the equivalence between (2.7) and (2.8) that,
	if $\mathbf{B}$ is of stable type $p$, then for every sequence $\{X, X_{n}; ~n \geq 1 \}$ 
	of i.i.d. {\bf B}-valued random variables, (2.3) holds if and only if
	\[
     \left \{
     \begin{array}{ll}
     \mbox{$\displaystyle \mathbb{E}X = 0 ~\mbox{and}~ \lim_{n \rightarrow \infty} 
     	n \mathbb{P}\left(\|X\| > n^{1/p} \right) = 0$} & \mbox{$\displaystyle \mbox{if}~ 1 < p < 2$,}\\
     &\\
    \mbox{$\displaystyle \lim_{n \rightarrow \infty} \left(\|\mathbb{E}XI(\|X\| \leq n) \| + 
    	n \mathbb{P}\left(\|X\| > n \right) \right) = 0$} & \mbox{$\displaystyle \mbox{if}~ p = 1$.}
    \end{array}
    \right.
	\]  
\end{remark}

\vskip 0.2cm

 It is well known that all real separable Hilbert spaces, $\mathcal{L}_{p}$ spaces ($p \geq 2$), and real separable finite-dimensional 
 Banach spaces are of stable type $p$ for all $0 < p \leq 2$. It now follows from Theorem 2.1 and Remarks 2.1 and 2.2 that the following 
 theorem provides large deviations for sums of heavy-tailed i.i.d. real-valued random variables.
 
 \vskip 0.2cm
 
 \begin{theorem}
  	Let $0 < p < 2$. Let $\{X, X_{n}; n \geq 1\}$ be a sequence of i.i.d. real-valued random variables such that
  	\begin{equation}
  	\left \{
  		\begin{array}{ll}
  		\mbox{$\displaystyle \lim_{n \rightarrow \infty} 
  			n \mathbb{P}\left(|X| > n^{1/p} \right) = 0$} & \mbox{$\displaystyle \mbox{if}~ 0 < p < 1$,}\\
  		&\\
  		\mbox{$\displaystyle \lim_{n \rightarrow \infty} \left(\mathbb{E}XI(|X|   \leq n) + 
  			n \mathbb{P}\left(|X| > n \right) \right) = 0$} & \mbox{$\displaystyle \mbox{if}~ p = 1$,}\\
  		&\\
  		\mbox{$\displaystyle \mathbb{E}X = 0 ~\mbox{and}~ \lim_{n \rightarrow \infty} 
  		n \mathbb{P}\left(|X| > n^{1/p} \right) = 0$} & \mbox{$\displaystyle \mbox{if}~ 1 < p < 2$.}\\
  		\end{array}
  		\right.
  \end{equation}
  	Then, for all $s > 0$,
  	\[
  	\limsup_{n \to \infty} \frac{1}{\log n} \log \mathbb{P}\left(\left|S_{n} \right| > s n^{1/p} \right)
  	= - \left(\bar{\beta} - p \right)/p
  	\]
  	and 
  	\[
  	\liminf_{n \to \infty} \frac{1}{\log n} \log \mathbb{P}\left(\left|S_{n} \right| > s n^{1/p} \right)
  	= -\left(\underline{\beta} - p \right)/p,
  	\]
  	where $\bar{\beta}$ and $\underline{\beta}$ are defined by (2.1) and (2.2) respectively when $\|X\|$
  	is replaced by $|X|$. Hence, under condition (2.9), for all $s > 0$,
  	\begin{equation}
  	\lim_{n \to \infty} \frac{1}{\log n} \log \mathbb{P}\left(\left|S_{n} \right| > s n^{1/p} \right)
  	= - \left(\hat{\beta} - p \right)/p ~~\mbox{if and only if}~~\bar{\beta} = \underline{\beta} = \hat{\beta}.
  	\end{equation}
 \end{theorem}

\vskip 0.2cm

Using Theorem 2.2, we can establish large deviations for sums of heavy-tailed i.i.d. non-negative random variables.
 
 \vskip 0.2cm
 
 \begin{theorem}
 	 Let $\{X, X_{n}; n \geq 1\}$ be a sequence of non-negative i.i.d. random variables. Write
 	 \[
 	 \bar{\alpha} = - \limsup_{t \rightarrow \infty} \frac{1}{t} \log \mathbb{P}(\log X > t)
 	 ~~\mbox{and}~~
 	 \underline{\alpha} = - \liminf_{t \rightarrow \infty} \frac{1}{t} \log \mathbb{P}(\log X > t). 
 	 \]
 	 (i)~~If $0 < p < 1$ and
 	 \begin{equation}
 	 \lim_{n \rightarrow \infty} n \mathbb{P}\left(X > n^{1/p} \right) = 0,
 	 \end{equation}
 	 then, for all $s > 0$,
 	 \[
 	 \limsup_{n \rightarrow \infty} 
 	 \frac{1}{\log n} \log \mathbb{P}\left(S_{n} > s n^{1/p} \right) = - \left(\bar{\alpha} - p \right)/p
 	 \]
 	 and
 	 \[
 	 \liminf_{n \rightarrow \infty} 
 	 \frac{1}{\log n} \log \mathbb{P}\left(S_{n} > s n^{1/p} \right) = - \left(\underline{\alpha} - p \right)/p.
 	 \]
 	  Hence, under condition (2.11), for all $s > 0$,
 	  \[ 	  
 	  \lim_{n \to \infty} \frac{1}{\log n} \log \mathbb{P}\left(S_{n} > s n^{1/p} \right)
 	  = - \left(\hat{\alpha} - p\right)/p ~~\mbox{if and only if}~~\bar{\alpha} = \underline{\alpha} = \hat{\alpha}.
 	  \]
 	   (ii)~~If $1 \leq p < 2$,
 	   \begin{equation}
 	   \mathbb{E}X = \mu \in (0, \infty),~~\mbox{and}~~\lim_{n \rightarrow \infty} n \mathbb{P}\left(X > n^{1/p} \right) = 0,
 	   \end{equation}
 	   then, for all $s > 0$,
 	   \begin{equation}
 	   \limsup_{n \rightarrow \infty} \frac{1}{\log n} \log \mathbb{P}\left(S_{n} > n \mu + s n^{1/p} \right) = - \left(\bar{\alpha} - p \right)/p
 	   \end{equation}
 	   and
 	   \begin{equation}
 	   \liminf_{n \rightarrow \infty} \frac{1}{\log n} \log \mathbb{P}\left(S_{n} > n \mu + s n^{1/p} \right) = - \left(\underline{\alpha} - p \right)/p.
 	   \end{equation}
 	   Hence, under condition (2.12), for all $s > 0$,
 	   \[ 	  \lim_{n \to \infty} \frac{1}{\log n} \log \mathbb{P}\left(S_{n} > n \mu + s n^{1/p} \right)
 	   = - \left(\hat{\alpha} - p \right)/p ~~\mbox{if and only if}~~\bar{\alpha} = \underline{\alpha} = \hat{\alpha}.
 	   \]
 \end{theorem}
 
 \vskip 0.2cm
 
 \begin{remark}
 	Theorem 2.3 (ii) is new. 

 	Theorem 2.3 (i) improves Theorem 1.2 (i.e., Theorems 2.1 and 2.2 of Hu and Nyrhinen (2004)). Note that, by Remark 3.1, 
 	\[
 	\bar{\alpha} = - \limsup_{t \rightarrow \infty} \frac{1}{t} \log \mathbb{P}(\log X > t) 
 	= \sup \left \{r \geq 0;~\lim_{n \rightarrow \infty} n \mathbb{P}\left(X > n^{1/r} \right) = 0 \right \}.
 	\]
 	If $\bar{\alpha} > 0$, then for every $x > \max\{1, 1/\bar{\alpha}\}$, we have $0 < 1/x < \min\{1, \bar{\alpha}\}$
 	and hence (2.11) holds with $p = 1/x$. Hence, by Theorem 2.3 (i), Theorem 1.2 (i) and Theorem 1.2 (iii) follow 
 	if $0 < \bar{\alpha} < \infty$ and $\bar{\alpha} = \infty$ respectively. 
 	
 	Obviously, the case where $\bar{\alpha} \in (0, \infty)$ and $\underline{\alpha} = \infty$ is not covered by 
 	Theorem 1.2. However, for this case, by Theorem 2.3 (i), for all $x > \max\{1, 1/\bar{\alpha}\}$,
 	\[
 	\limsup_{n \rightarrow \infty} 
 	\frac{1}{\log n} \log \mathbb{P}\left(S_{n} > s n^{x} \right) = \limsup_{n \rightarrow \infty} 
 	\frac{1}{\log n} \log \mathbb{P}\left(S_{n} > s n^{1/p} \right) =- \left(\bar{\alpha} - p \right)/p
 	= - \bar{\alpha}x + 1
 	\]
 	and
 	\[
 	\liminf_{n \rightarrow \infty} 
 	\frac{1}{\log n} \log \mathbb{P}\left(S_{n} > s n^{x} \right) = \liminf_{n \rightarrow \infty} 
 	\frac{1}{\log n} \log \mathbb{P}\left(S_{n} > s n^{1/p} \right) = - \infty.
 	\]
 	
 	The assertion of Theorem 1.2 (ii) is simple. In fact, since $\{X, X_{n}; n \geq 1\}$ is a sequence of 
 	non-negative i.i.d. random variables, for all $s > 0$,
 	\[
 	\mathbb{P}\left(S_{n} > s n^{x} \right) \geq \mathbb{P}\left(\max_{1 \leq k \leq n} X_{k} > s n^{x} \right)
 	= 1 - \left(1 - \mathbb{P}\left(X > s n^{x} \right) \right)^{n}.
 	\]
 	If $\bar{\alpha} = 0$, then, by Remark 3.1, for all $x > 0$ and all $s > 0$,
 	\[
 	\limsup_{n \rightarrow \infty} n \mathbb{P}\left(X > s n^{x} \right) = \infty
 	\]
 	which ensures
 	\[
 	\liminf_{n \rightarrow \infty} \left( 1 - \mathbb{P}\left(X > s n^{x} \right) \right)^{n} = 0, 
 	\] 
 	and hence, 
 	\[
 	\limsup_{n \rightarrow \infty} \mathbb{P} \left(S_{n} > s n^{x} \right) = 1,
 	\]
 	i.e., Theorem 1.2 (ii) follows.
 	\end{remark}
 
 \section{Preliminary lemmas}	
 	
 To prove our main results, we use the following preliminary lemmas.
 	
 \vskip 0.2cm
 	
 \begin{lemma}
 	Let $X$ be a $\mathbf{B}$-valued random variable. Then 
 	\begin{equation}
 	\bar{\beta} = - \limsup_{t \rightarrow \infty} \frac{1}{t} \log \mathbb{P}(\log \|X\| > t) 
 	= \sup \left \{r \geq 0;~\lim_{t \rightarrow \infty} t^{r} \mathbb{P}(\|X\| > t) = 0 \right \},
 	\end{equation}
 	\begin{equation}
 	\underline{\beta} = - \liminf_{t \rightarrow \infty} \frac{1}{t} \log \mathbb{P}(\log \|X\| > t) 
 	 = \sup \left \{r \geq 0;~\liminf_{t \rightarrow \infty} t^{r} \mathbb{P}(\|X\| > t) = 0 \right \}.
 	 	\end{equation}
 \end{lemma}
 
 \vskip 0.2cm
 
 \noindent {\it Proof}~~We will only give the proof of (3.1) since a proof of (3.2) can be culled from 
 the proof of (3.1) with obvious modifications.
 
 Replacing $t$ by $\log t$, we have 
 \[
 \bar{\beta} = - \limsup_{t \rightarrow \infty} \frac{1}{t} \log \mathbb{P}(\log \|X\| > t) 
 = - \limsup_{t \rightarrow \infty} \frac{1}{\log t} \log \mathbb{P}( \|X\| > t). 
 \]
If $0 \leq \bar{\beta} < \infty$, then, for any given $\epsilon > 0$, 
\begin{equation}
\mathbb{P}(\|X\| > t) \leq t^{-\bar{\beta} + \epsilon} ~~\mbox{sufficiently large}~ t
\end{equation}
and there exists a sequence $\left \{t_{n, \epsilon};~ n \geq 1 \right \}$ of  increasing positive 
numbers such that
\begin{equation}
\lim_{n \rightarrow \infty} t_{n, \epsilon} = \infty
~~\mbox{and}~~\mathbb{P}\left(\|X\| > t_{n, \epsilon} \right) \geq t_{n, \epsilon}^{-\bar{\beta} - \epsilon}, ~~n \geq 1.
\end{equation}
From (3.3) and (3.4), we have
\[
\bar{\beta} - \epsilon \leq \sup \left \{r \geq 0;~\lim_{t \rightarrow \infty} t^{r} \mathbb{P}(\|X\| > t) = 0 \right \} \leq \bar{\beta} + \epsilon,
\]
and hence, letting $\epsilon \searrow 0$, (3.1) holds if $0 \leq \bar{\beta} < \infty$.

If $\bar{\beta} = \infty$, then, for any given $M > 0$,
\[
\mathbb{P}(\|X\| > t) \leq t^{-M} ~~\mbox{sufficiently large}~ t
\]
which ensures 
\[
\sup \left \{r \geq 0;~\lim_{t \rightarrow \infty} t^{r} \mathbb{P}(\|X\| > t) = 0 \right \} \geq M.
\]
Hence, letting $M \rightarrow \infty$, (3.1) follows if $\bar{\beta} = \infty$. ~~$\Box$

\vskip 0.2cm

\begin{remark}
	Let $X$ be a $\mathbf{B}$-valued random variable. Let $\bar{\beta}$ and $\underline{\beta}$ 
	be defined in Lemma 3.1. By Lemma 3.1, it is easy to see that, for all $s > 0$,
	\[
	\bar{\beta} =  \sup \left \{r \geq 0;~\lim_{t \rightarrow \infty} t^{r} \mathbb{P}(\|X\| > t) = 0 \right \}
	= \sup \left \{r \geq 0;~\lim_{t \rightarrow \infty} t^{r} \mathbb{P}(\|X\| > st) = 0 \right \}
	\]
and
	\[
	\underline{\beta} 
	=  \sup \left \{r \geq 0;~\liminf_{t \rightarrow \infty} t^{r} \mathbb{P}(\|X\| > t) = 0 \right \}
	= \sup \left \{r \geq 0;~\liminf_{t \rightarrow \infty} t^{r} \mathbb{P}(\|X\| > st ) = 0 \right \},
	\]
and hence,
	\[
	\bar{\beta} =  \sup \left \{r \geq 0;~\lim_{t \rightarrow \infty} t \mathbb{P}\left(\|X\| > st^{1/r} \right) = 0 \right \}
	= \sup \left \{r \geq 0;~\lim_{n \rightarrow \infty} n \mathbb{P}\left(\|X\| > sn^{1/r} \right) = 0 \right \}
	\]
	and
	\[
	\underline{\beta} 
	=  \sup \left \{r \geq 0;~\liminf_{t \rightarrow \infty} t \mathbb{P}\left(\|X\| > st^{1/r} \right) = 0 \right \}
	= \sup \left \{r \geq 0;~\liminf_{n \rightarrow \infty} n \mathbb{P}\left(\|X\| > sn^{1/r} \right) = 0 \right \}.
	\]
\end{remark}

\vskip 0.2cm

\begin{remark}
	Let $X$ be a $\mathbf{B}$-valued random variable. Let $X^{\prime}$ be an independent copy of $X$. Write
	$t_{0} = \inf\left\{t \geq 0;~ \mathbb{P}(\|X\| > t ) \leq 1/2 \right \}$. Then $t_{0} < \infty$. Hence, for $t > t_{0}$, as
	\[
	\left\{\|X\| > 2t, \|X^{\prime}\| \leq t \right\} \subseteq \left\{\|X - X^{\prime} \| > t \right \} \subseteq \left\{\|X\| > t/2 
	\right\},
	\]
	we have
	\[
	\mathbb{P}(\|X\| > 2t) \leq 2 \mathbb{P}\left(\|X - X^{\prime} \| > t \right).
	\]
	Then, for $r \geq 0$,
	\[
	\lim_{t \rightarrow \infty} t^{r} \mathbb{P}(\|X\| > t) = 0
	~~\mbox{if and only if}~~\lim_{t \rightarrow \infty} t^{r} \mathbb{P}\left(\|X - X^{\prime}\| > t\right) = 0
	\]
	and
	\[
	\liminf_{t \rightarrow \infty} t^{r} \mathbb{P}(\|X\| > t) = 0
	~~\mbox{if and only if}~~\liminf_{t \rightarrow \infty} t^{r} \mathbb{P}\left(\|X - X^{\prime}\| > t\right) = 0.
	\]
	Thus, by Lemma 3.1, we have
	\[
	 \bar{\beta} = -\limsup_{t \rightarrow \infty} \frac{1}{t} \log \mathbb{P}(\log \|X\| > t) =
	-\limsup_{t \rightarrow \infty} \frac{1}{t} \log \mathbb{P}(\log \|X - X^{\prime}\| > t) 
	\] 
	and
	\[
	 \underline{\beta} = -\liminf_{t \rightarrow \infty} \frac{1}{t} \log \mathbb{P}(\log \|X\| > t) =
	 -\liminf_{t \rightarrow \infty} \frac{1}{t} \log \mathbb{P}(\log \|X - X^{\prime}\| > t).
	 \] 
\end{remark}

\vskip 0.2cm

The following lemma, which provides a symmetrization procedure for the large deviation for $\mathbf{B}$-valued 
random variables, may be of independent interest.

\vskip 0.2cm

\begin{lemma}
	Let $\{Y_{n};~ n \geq 1 \}$ be a sequence of $\mathbf{B}$-valued random variables such that
	\begin{equation}
	Y_{n} \rightarrow_{\mathbb{P}} 0.
	\end{equation}
	Let $\{Y_{n}^{\prime};~ n \geq 1 \}$ be an independent copy of $\{Y_{n}; n \geq 1\}$. 
	Write $\hat{Y}_{n} = Y_{n} - Y^{\prime}_{n}$, $n \geq 1$. Let $\{a_{n};~n \geq 1 \}$ be a sequence of positive 
	numbers such that $\lim_{n \rightarrow \infty} a_{n} = \infty$. If there exist two constants 
	$0 \leq \bar{\gamma} \leq \underline{\gamma} \leq \infty$ such that, for all $s > 0$,
	\begin{equation}
	\limsup_{n \rightarrow \infty} \frac{1}{a_{n}} \log \mathbb{P}\left(\|\hat{Y}_{n} \| > s \right) = - \bar{\gamma}
	~~\mbox{and}~~
	\liminf_{n \rightarrow \infty} \frac{1}{a_{n}} \log \mathbb{P}\left(\|\hat{Y}_{n}\| > s \right) = - \underline{\gamma},
	\end{equation}
	then, for all $s > 0$,
	\begin{equation}
	\limsup_{n \rightarrow \infty} \frac{1}{a_{n}} \log \mathbb{P}\left(\|Y_{n} \| > s \right) = - \bar{\gamma}
	~~\mbox{and}~~
	\liminf_{n \rightarrow \infty} \frac{1}{a_{n}} \log \mathbb{P}\left(\|Y_{n}\| > s \right) = - \underline{\gamma}.
	\end{equation}
\end{lemma}

\vskip 0.2cm

\noindent {\it Proof}~~For any given $s > 0$, it follows from (3.5) that there exists a positive integer $n_{0}$ (which depends
on $s$ only) such that, for all $n \geq n_{0}$, 
\[
\mathbb{P}\left(\|Y_{n}\| > s/2 \right) \leq 1/2.
\]
Then, for all $n \geq n_{0}$,
\[
\left\{\|Y_{n}\| > s, \|Y^{\prime}_{n}\| \leq s/2 \right\} \subseteq  \left\{\|\hat{Y}_{n} \| > s/2 \right \} 
~~\mbox{and}~~
\left\{\|\hat{Y}_{n} \| > 2s \right \} \subseteq \left\{\|Y_{n}\| > s \right\} \cup \left\{\|Y^{\prime}_{n} \| > s \right\}.
\]
Then, for all $n \geq n_{0}$,
\[
\frac{1}{2} \mathbb{P}\left( \|\hat{Y}_{n} \| > 2s \right) \leq \mathbb{P} 
\left( \|Y_{n} \| > s \right) \leq 2 \mathbb{P}\left( \|\hat{Y}_{n} \| > s/2 \right),
\]
and hence, 
\begin{equation}
\frac{1}{a_{n}} \log \left(\frac{1}{2} \mathbb{P}\left( \|\hat{Y}_{n} \| > 2s \right) \right) 
\leq \frac{1}{a_{n}} \log \mathbb{P}\left( \|Y_{n} \| > s \right) 
\leq \frac{1}{a_{n}} \log \left(2 \mathbb{P}\left( \|\hat{Y}_{n} \| > s/2 \right) \right).
\end{equation}
Since $\lim_{n \rightarrow \infty} a_{n} = \infty$, (3.7) follows from (3.8) and (3.6).~~$\Box$

\vskip 0.2cm

The first part of the following lemma is one of L\'{e}vy's inequalities
(see, e.g., Ledoux and Talagrand (1991, p.47)), the second part is a
version of the Hoffmann-J{\o}rgensen inequality due to Li, Rao,
Jiang, and Wang (1995, Lemma 2.2), and the third part is the second part of
Proposition 6.8 of Ledoux and Talagrand (1991, p.156).   

\begin{lemma}
	Let $\{V_{k};~1\leq k\leq n\}$ be a finite sequence of independent symmetric 
	$\mathbf{B}$-valued random variables, and set $T_{n}=V_{1}+...+V_{n}$. We have\\
	(i)~~For any $t>0$,
	\[
	\mathbb{P}\left( \max_{1\leq k\leq n}\|V_{k}\|>t\right) \leq 2 \mathbb{P}\left( \|T_{n}\|>t\right);
	\]
	(ii)~~For each integer $j\geq 1$, there exist positive numbers $C_{j}$ and $D_{j}$, depending only on $j$,
	such that
	\[
	\mathbb{P}\left(\|T_{n}\|>2jt \right) \leq  C_{j} \mathbb{P}\left( \max_{1\leq k\leq n}\|V_{k}\| > t\right)
	+ D_{j}\left( \mathbb{P} (\|T_{n}\|>t) \right) ^{j};
	\]
	(iii)~If, for some $0 < r < \infty$, $\mathbb{E}\left \|V_{k} \right \|^{r} < \infty$, $1 \leq k \leq n$, 
	then
	\[
	\mathbb{E}\left \|T_{n} \right \|^{r} 
	\leq 2\cdot 3^{r} \mathbb{E}\left(\max_{1 \leq k \leq n}\|V_{k} \|^{r} \right) + 2\left(3t_{n}\right)^{r},
	\]
	where $\displaystyle t_{n} = \inf\left\{t > 0;~\mathbb{P}\left(\|T_{n} \|  > t \right) \leq \left(8 \cdot 3^{r} \right)^{-1} \right\}$.
\end{lemma}

\vskip 0.2cm

The following lemma improves Corollary 2.3 of Hu and Nyrhinen (2004).

\vskip 0.2cm

\begin{lemma}
	Let $\{X, X_{n}; n \geq 1\}$ be a sequence of i.i.d. $\mathbf{B}$-valued random variables. Let $p > 0$. Then, for any $s > 0$,
	\begin{equation}
	\limsup_{n \rightarrow \infty} \frac{1}{\log n} \log \mathbb{P}\left(\max_{1 \leq k \leq n} \|X_{k}\| > sn^{1/p} \right)
	= 0 \wedge \left( - (\bar{\beta} - p)/p \right)
	\end{equation}
	and
	\begin{equation}
	\liminf_{n \rightarrow \infty} \frac{1}{\log n} \log \mathbb{P}\left(\max_{1 \leq k \leq n} \|X_{k}\| > sn^{1/p} \right)
	= 0 \wedge \left( - (\underline{\beta} - p)/p \right).
	\end{equation}
\end{lemma}

\vskip 0.2cm

\noindent {\it Proof}~~For $n \geq 1 $, let $h_{n}(x) = 1 - (1-x)^{n}$, $0 \leq x \leq 1$. Clearly, $h_{n}(x)$ 
is an increasing function on $[0, 1]$. Then, for $1/n \leq x \leq 1$, we have
\[
\frac{1}{2} \leq 1 - \left(1 - \frac{1}{n} \right)^{n} = h_{n}(1/n) \leq 1 - (1 - x)^{n} \leq h_{n}(1) = 1.
\]
For $0 \leq x \leq 1/n$, using the second-degree Taylor polynomial of $h_{n}(x)$ at $x = 0$ with remainder, we have
\[
\frac{nx}{2} \leq nx - \frac{(nx)^{2}}{2} \leq nx - \frac{n(n-1)x^{2}}{2} \leq 1 - (1-x)^{n} \leq nx.
\]
Hence, for all $0 \leq x \leq 1$, we have
\begin{equation}
\frac{1 \wedge (nx)}{2} \leq 1 - (1 - x)^{n} \leq 1 \wedge (nx).
\end{equation}
Note that
\[
\mathbb{P}\left(\max_{1 \leq k \leq n} \|X_{k}\| > sn^{1/p} \right)
= 1 - \left(1 - \mathbb{P}\left(\|X\| > sn^{1/p} \right)\right)^{n}, ~n \geq 1.
\]
Thus it follows from (3.11) that
\[
\frac{1 \wedge \left(n \mathbb{P}\left(\|X\| > sn^{1/p} \right) \right)}{2}
\leq \mathbb{P}\left(\max_{1 \leq k \leq n} \|X_{k}\| > sn^{1/p} \right)
\leq 1 \wedge \left(n \mathbb{P}\left(\|X\| > sn^{1/p} \right)\right), ~ n \geq 1,
\]
and hence,
\[
\begin{array}{ll}
& \mbox{$\displaystyle 
\limsup_{n \rightarrow \infty} \frac{1}{\log n} 
\log \mathbb{P}\left(\max_{1 \leq k \leq n} \|X_{k}\| > sn^{1/p} \right)$}\\
&\\
& \mbox{$\displaystyle =  
\limsup_{n \rightarrow \infty} \frac{1}{\log n} 
\log \left(1 \wedge \left(n \mathbb{P}\left(\|X\| > sn^{1/p} \right)\right) \right)$}\\
&\\
& \mbox{$\displaystyle = 
	0 \wedge \left(\limsup_{n \rightarrow \infty} 
	\frac{1}{\log n} \log \left(n \mathbb{P}\left(\|X\| > sn^{1/p} \right) \right) \right)$}\\
&\\
& \mbox{$\displaystyle = 
	0 \wedge \left(1 + \limsup_{n \rightarrow \infty} 
	\frac{1}{\log n} \log \mathbb{P}\left(\|X\| > sn^{1/p} \right) \right)$}\\
&\\
& \mbox{$\displaystyle = 
	0 \wedge \left(1 + \limsup_{n \rightarrow \infty} 
	\left(\frac{\log \left(sn^{1/p} \right)}{\log n} 
	\times \frac{1}{\log\left(sn^{1/p} \right)} \log \mathbb{P}\left(\|X\| > sn^{1/p} \right) \right) \right)$}\\
&\\
& \mbox{$\displaystyle = 
	0 \wedge \left(1 + \frac{1}{p}\limsup_{t \rightarrow \infty} \frac{1}{\log t} \log \mathbb{P}(\|X\| > t) \right)$}\\
&\\
&  \mbox{$\displaystyle = 
0 \wedge \left(1 + \frac{1}{p}\left(-\bar{\beta}\right) \right)$}\\
&\\
&  \mbox{$\displaystyle = 0 \wedge \left( - (\bar{\beta} - p)/p \right)$},
\end{array}
\]
proving (3.9). Similary (3.10) follows when $\displaystyle \limsup_{n\rightarrow \infty}$ 
is replaced by $\displaystyle \liminf_{n\rightarrow \infty}$.  $\Box$

\vskip 0.2cm

\begin{lemma}
	Let $0 < p < 2$. Let $\{X, X_{n}; n \geq 1\}$ be a sequence of i.i.d. $\mathbf{B}$-valued random variables 
	such that (2.3) holds. Then
	\begin{equation}
		\lim_{t \rightarrow \infty} t^{p} \mathbb{P}(\|X\| > t) = 0
		~~\mbox{and}~~p \leq \bar{\beta}.
	\end{equation}
\end{lemma}

\vskip 0.2cm

\noindent {\it Proof}~~Let $\{X^{\prime}, X_{n}^{\prime};~ n \geq 1 \}$ 
be an independent copy of $\{X, X_{n}; n \geq 1\}$. Let $\hat{X} = X - X^{\prime}$ 
and $\hat{X}_{n} = X_{n} - X^{\prime}_{n}$, $n \geq 1$. Then 
$\{\hat{X}, \hat{X}_{n}; n \geq 1\}$ is a sequence of i.i.d. 
symmetric $\mathbf{B}$-valued random variables. Note that (2.3) implies
\[
\frac{\hat{S}_{n}}{n^{1/p}} \rightarrow_{\mathbb{P}} 0,
\]
where $\hat{S}_{n} = \sum_{k=1}^{n} \hat{X}_{k}$, $n \geq 1$. Thus, applying Lemma 3.3 (i), we have
\[
\mathbb{P}\left(\max_{1 \leq k \leq n} \|\hat{X}_{k} \| > n^{1/p} \right) 
\leq 2 \mathbb{P} \left(\|\hat{S}_{n} \| > n^{1/p} \right) \rightarrow 0 ~~\mbox{as}~~n \rightarrow \infty,
\]
and hence,
\[
\mathbb{P}\left(\max_{1 \leq k \leq n} \|\hat{X}_{k} \| > n^{1/p} \right) 
= 1 - \left(1 - \mathbb{P}\left(\|\hat{X}\| > n^{1/p} \right) \right)^{n}
\rightarrow 0~~\mbox{as}~~n \rightarrow \infty
\]
which is equivalent to
\[
n \mathbb{P}\left(\|\hat{X}\| > n^{1/p} \right) \rightarrow 0~~\mbox{as}~~n \rightarrow \infty.
\]
Thus, by Remarks 3.1 and 3.2, (3.12) follows.~~$\Box$

\vskip 0.2cm

The following lemma is one of the de Acosta (1981) inequalities.

\vskip 0.2cm

\begin{lemma}
	For every $1 \leq r \leq 2$ there exists a positive constant $C(r)$ such that for any 
	separable Banach space $\mathbf{B}$ and any finite sequence $\{V_{k};~ 1 \leq k \leq n \}$ 
	of independent $\mathbf{B}$-valued random variables with 
	$\mathbb{E}\|V_{k}\|^{r} < \infty$, $1 \leq k \leq n$, and $T_{n}=V_{1}+...+V_{n}$,
	\[
	\mathbb{E} \left|\| T_{n} \| - \mathbb{E} \|T_{n} \| \right|^{r} \leq C(r) \sum_{k=1}^{n} \mathbb{E} \|V_{k}\|^{r};
	\]
	if $r = 2$ then it is possible to take $C(2) = 4$. 
\end{lemma}

\vskip 0.2cm

Let $m(Y)$ denote a median for a real-valued random variable $Y$. The following lemma is used to prove Theorem 2.3. 
The first part was obtained by Li and Rosalsky (2004, Lemma 3.1) and the second part was established by Petrov (1975) 
(also see Petrov (1995, Theorem 2.1)).

\vskip 0.2cm

\begin{lemma}
	Let $\{V_{k};~1 \leq k \leq n \}$ be a finite sequence of independent real-valued random variables and set
	$T_{0} = 0$ and $T_{k} = V_{1} + \cdots + V_{k}$, $1 \leq k \leq n$. Then, for every real $t$,
	\begin{equation}
	\mathbb{P}\left(\max_{1 \leq k \leq n} \left(V_{k} + m\left(T_{k-1}\right) \right) > t \right) 
	\leq 2 \mathbb{P}\left(\max_{1 \leq k \leq n} T_{k} > t \right),
	\end{equation}
	\begin{equation}
	\mathbb{P}\left(\max_{1 \leq k \leq n} \left(T_{k} + m\left(T_{n} - T_{k}\right) \right) > t \right) 
	\leq 2 \mathbb{P}\left(T_{n} > t \right).
	\end{equation}
	
\end{lemma}

\section{Proofs of the main results}

With the preliminaries accounted for, Theorem 2.1 may be proved.

\vskip 0.3cm

{\it Proof of Theorem 2.1}~~Applying Remark 3.2 and Lemma 3.2, without loss of generality, we can assume that 
 $\{X, X_{n}; n \geq 1\}$ is a sequence of i.i.d. symmetric $\mathbf{B}$-valued random variables. Thus, by Lemma 3.3 (i), for all $s > 0$, 
 \[
 \mathbb{P}\left(\max_{1 \leq k \leq n} \|X_{k}\| > sn^{1/p} \right) \leq 2 \mathbb{P}\left(\|S_{n}\| > sn^{1/p} \right), ~n \geq 1.
 \]
 Hence, by Lemma 3.4, we have
 \[
 \limsup_{n \to \infty} \frac{1}{\log n} \log \mathbb{P}\left(\left\|S_{n} \right\| > s n^{1/p} \right) 
 \geq \limsup_{n \to \infty} \frac{1}{\log n} \log \mathbb{P}\left( \max_{1 \leq k \leq n}\left\|X_{k} \right\| > s n^{1/p} \right)
 = 0 \wedge \left( - (\bar{\beta} - p)/p \right)
 \]
 and
 \[
 \liminf_{n \to \infty} \frac{1}{\log n} \log \mathbb{P}\left(\left\|S_{n} \right\| > s n^{1/p} \right) 
 \geq \liminf_{n \to \infty} \frac{1}{\log n} \log \mathbb{P}\left( \max_{1 \leq k \leq n}\left\|X_{k} \right\| > s n^{1/p} \right)
 = 0 \wedge \left( - (\underline{\beta} - p)/p \right).
 \] 
 By Lemma 3.5, it follows from (2.3) that $p \leq \bar{\beta} \leq \underline{\beta}$. Thus,
 \begin{equation}
 \limsup_{n \to \infty} \frac{1}{\log n} \log \mathbb{P}\left(\left\|S_{n} \right\| > s n^{1/p} \right) 
 \geq - (\bar{\beta} - p)/p
 \end{equation}
 and 
  \begin{equation}
  \liminf_{n \to \infty} \frac{1}{\log n} \log \mathbb{P}\left(\left\|S_{n} \right\| > s n^{1/p} \right) 
  \geq - (\underline{\beta} - p)/p.
  \end{equation}
 
 For given $p < q < \infty$, write
 \[
 V_{n,k} = X_{k}I\left(\|X_{k}\| \leq n^{1/q} \right), ~1 \leq k \leq n, ~T_{n} = \frac{\sum_{k=1}^{n}V_{n,k}}{n^{1/p}},~~n \geq 1.
 \]
 Since, for each $n \geq 1$, $X_{1}, ..., X_{n}$ are symmetric i.i.d. $\mathbf{B}$-valued random variables, it is easy to see that 
 \[
 \left\{X_{k} I\left(|X_{k}| \leq n^{1/q} \right) -  X_{k} I\left(|X_{k}| > n^{1/q} \right);~1 \leq k \leq n \right \}
 ~~\mbox{and}~~\left\{X_{k};~ 1 \leq k \leq n \right \}
 \]
 are identically distributed finite sequence of symmetric i.i.d. $\mathbf{B}$-valued random variables. Note that, for each $n \geq 1$,
 \[
 T_{n} = \frac{1}{2} 
 \left(\frac{\sum_{k=1}^{n} \left(X_{k} I\left(\|X_{k}\| \leq n^{1/q} \right) -  X_{k} I\left(\|X_{k}\| > n^{1/q} \right) \right)}{n^{1/p}} 
 + \frac{S_{n}}{n^{1/p}} \right).
 \] 
 Thus it follows from (2.3) that, for any $\epsilon > 0$,
 \[
 \begin{array}{lll}
\mbox{$\displaystyle  
	\mathbb{P} \left(\left\|T_{n} \right \| > \epsilon \right)$}
& \leq & \mbox{$\displaystyle 
  \mathbb{P} \left(\left\|\frac{ \sum_{k=1}^{n} \left(X_{k} I\left(\|X_{k}\| \leq n^{1/q} \right) 
		-  X_{k} I\left(\|X_{k}\| > n^{1/q} \right) \right)}{n^{1/p}} \right \| > \epsilon \right) 
 + \mathbb{P} \left(\left\|\frac{S_{n}}{n^{1/p}} \right \| > \epsilon \right)$}\\ 
&&\\
& = & \mbox{$\displaystyle 
2 \mathbb{P} \left(\left\|\frac{S_{n}}{n^{1/p}} \right \| > \epsilon \right)$}\\ 
&&\\
& \rightarrow & \mbox{$\displaystyle 0~~\mbox{as}~~n \rightarrow \infty$}
\end{array}
\]
and hence,
\[
t_{n} = \inf\left\{t > 0;~ \mathbb{P} \left(\left\|T_{n} \right \| > t \right) \leq 1/24 \right \} 
\rightarrow 0~~\mbox{as}~~n \rightarrow \infty.
\]
By Lemma 3.3 (iii) with $r = 1$, we have
\begin{equation}
\begin{array}{lll}
\mbox{$\displaystyle
\mathbb{E}\left\|T_{n} \right\| $}
& \leq &
\mbox{$\displaystyle
	6 \mathbb{E} \left(\max_{1 \leq k \leq n} \left \|\frac{X_{k}I\left(\|X_{k}\|\leq n^{1/q} \right)}{n^{1/p}} \right \| \right) 
	+ 6t_{n}$}\\
&&\\
& \leq &
\mbox{$\displaystyle
	6\left(n^{-(q - p)/(pq)} + t_{n} \right)$}\\
&&\\
& \rightarrow &
\mbox{$\displaystyle 
	0~~\mbox{as}~~n \rightarrow \infty$.}
\end{array}
\end{equation}
Let $\delta = \frac{p}{2} \wedge \frac{(2-p)(q - p)}{2p} > 0$. Then 
\[
\tau = \frac{2}{p} -1 - \frac{2 - p + \delta}{q} = \frac{2q - pq - 2p + p^{2} - \delta p}{pq} 
= \frac{(q - p)(2 - p) - \delta p}{pq} \geq \frac{(q - p)(2 - p)}{2pq} > 0.
\]
By Lemma 3.5, (2.3) implies that 
\[
\mathbb{E}\|X\|^{p - \delta} < \infty.
\]
For given $s > 0$, it follows from (4.3), Markov's inequality, and Lemma 3.6 (with $r = 2$) that, for all sufficiently large $n$,
\begin{equation}
\begin{array}{lll}
\mbox{$\displaystyle 
\mathbb{P}\left(\|T_{n} \| > s \right)$}
& \leq &
\mbox{$\displaystyle 
	\mathbb{P}\left(\left|\|T_{n} \| - \mathbb{E}\|T_{n}\| \right| > \frac{s}{2} \right)$}\\
&&\\
& \leq &
\mbox{$\displaystyle \frac{\mathbb{E}\left(\|T_{n} \| - \mathbb{E}\|T_{n}\| \right)^{2}}{(s/2)^{2}}$}\\
&&\\
& \leq &
\mbox{$\displaystyle \frac{16}{s^{2}} \sum_{k=1}^{n}\frac{\mathbb{E}\left(\|V_{n,k}\|^{2}\right)}{n^{2/p}}$}\\
&&\\
&=& 
\mbox{$\displaystyle \left(\frac{16}{s^{2}}\right) \frac{n \mathbb{E}\left(\|X\|^{2}I\left(\|X\| \leq n^{1/q} \right) \right)}{n^{2/p}}$}\\
&&\\
&\leq & 
\mbox{$\displaystyle 
	\left(\frac{16}{s^{2}}\right) \frac{n \mathbb{E}\left(\|X\|^{p-\delta}\left(n^{1/q}\right)^{2-p+\delta} \right)}{n^{2/p}}$}\\
&&\\
& = & 
\mbox{$\displaystyle 
	\left(\frac{16\mathbb{E}\|X\|^{p-\delta}}{s^{2}}\right) n^{- \left(\frac{2}{p} - 1 - \frac{2 - p + \delta}{q}\right)}$}\\
&&\\
& = & 
\mbox{$\displaystyle \mbox{\it O}\left(n^{-\tau}\right)$}.
\end{array}
\end{equation}
For given $M > 0$, let $j = [M/\tau] + 1$ and $\displaystyle t = \frac{s}{2j}$. Then, it follows from Lemma 3.3 (ii) 
and (4.4) that, for all sufficiently large $n$,
\[
\begin{array}{lll}
\mbox{$\displaystyle 
	\mathbb{P}\left(\|T_{n} \| > s \right)$}
& = &
\mbox{$\displaystyle 
	\mathbb{P}\left(\|T_{n} \| > 2jt \right)$}\\
&&\\
&\leq&
\mbox{$\displaystyle 
	 C_{j} \mathbb{P}\left( \max_{1\leq k\leq n}\|V_{n,k}\| > t\right)
	 + D_{j}\left( \mathbb{P} (\|T_{n}\|>t) \right) ^{j}$}\\
&&\\
&\leq&
\mbox{$\displaystyle 
	C_{j} \mathbb{P}\left(n^{-\frac{q-p}{pq}} > t \right) + \mbox{\it O}\left(n^{-j\tau} \right) $}\\
&=&
\mbox{$\displaystyle \mbox{\it O}\left(n^{-M} \right) $},
\end{array}
\]
and hence,
\[
\limsup_{n \rightarrow \infty} \frac{1}{\log n} \log \mathbb{P}\left(\|T_{n} \| > s \right) \leq -M.
\]
Letting $M \rightarrow \infty$, we have, for all $s > 0$,
\begin{equation}
\lim_{n \rightarrow \infty} \frac{1}{\log n} \log \mathbb{P}\left(\|T_{n} \| > s \right) = - \infty.
\end{equation}
Since, for all $s > 0$,
\[
\begin{array}{lll}
\mbox{$\displaystyle 
\mathbb{P}\left(\|S_{n}\| > sn^{1/p} \right)$}
&\leq&
\mbox{$\displaystyle 
 \mathbb{P}\left(\|T_{n}\| > s \right) + \mathbb{P} \left(\frac{S_{n}}{n^{1/p}} - T_{n} \neq 0 \right)$}\\
&&\\
&=&
\mbox{$\displaystyle 
	\mathbb{P}\left(\|T_{n}\| > s \right) + \mathbb{P} \left(\frac{\sum_{k=1}^{n}X_{k}I\left(\|X_{k}\| > n^{1/q} \right)}{n^{1/p}} \neq 0 \right)$}\\
&&\\
&\leq&
\mbox{$\displaystyle 
	\mathbb{P}\left(\|T_{n}\| > s \right) + n \mathbb{P}\left(\|X\| > n^{1/q} \right)$}\\
&&\\
&\leq&
\mbox{$\displaystyle 
   \left(2\mathbb{P}\left(\|T_{n}\| > s \right) \right) \vee \left(2 n \mathbb{P}\left(\|X\| > n^{1/q} \right)\right),~~n \geq 1,$}\\
\end{array}
\]
we have
\[
\log \mathbb{P}\left(\|S_{n}\| > sn^{1/p} \right) \leq \log 2 
+ \left(\log \mathbb{P}\left(\|T_{n}\| > s \right) \right) 
\vee  \left(\log \left( n \mathbb{P}\left(\|X\| > n^{1/q} \right)\right) \right),~~n \geq 1.
\]
Note that
\[
\begin{array}{lll}
\mbox{$\displaystyle
\limsup_{n \rightarrow \infty} \frac{1}{\log n} \log \left(n \mathbb{P}\left(\|X\| > n^{1/q} \right)\right)$}
&=&
\mbox{$\displaystyle
	1 + \limsup_{t \rightarrow \infty} \frac{1}{\log t} \log \left(\mathbb{P}\left(\|X\| > t^{1/q} \right)\right)$}\\
&&\\
&=&
\mbox{$\displaystyle
	1 + \limsup_{t \rightarrow \infty} \left(\frac{\log t^{1/q}}{\log t} \right) 
	\frac{1}{\log t^{1/q}} \log \left(\mathbb{P}\left(\|X\| > t^{1/q} \right)\right)$}\\
&&\\
&=&
\mbox{$\displaystyle
	 - \left(\bar{\beta} - q \right)/q$}
\end{array}
\]
and, similarly,
\[
\liminf_{n \rightarrow \infty} \frac{1}{\log n} \log \left(n \mathbb{P}\left(\|X\| > n^{1/q} \right)\right) 
=  - \left(\underline{\beta} - q \right)/q.
\]
Thus, it follows from (4.5) that, for all $s > 0$,
\[
\begin{array}{ll}
& \mbox{$\displaystyle
\limsup_{n \rightarrow \infty} \frac{1}{\log n} \log \mathbb{P}\left(\|S_{n}\| > s n^{1/p} \right)$}\\
&\\
& \mbox{$\displaystyle 
	\leq \limsup_{n \rightarrow \infty} \frac{1}{\log n} \left(\log 2 
	+ \left(\log \mathbb{P}\left(\|T_{n}\| > s \right) \right) 
	\vee  \left(\log \left( n \mathbb{P}\left(\|X\| > n^{1/q} \right)\right) \right)\right)$}\\
&\\
& \mbox{$\displaystyle 
	\leq \left(\lim_{n \rightarrow \infty} \frac{1}{\log n} \log \mathbb{P}\left(\|T_{n}\| > s \right) \right) \vee
	\left(\limsup_{n \rightarrow \infty} \frac{1}{\log n} \log \left(n \mathbb{P}\left(\|X\| > n^{1/q} \right) \right)\right)$}\\
&\\
& \mbox{$\displaystyle 
	= (-\infty) \vee \left((\bar{\beta} - q)/q \right)$}\\
&\\
& \mbox{$\displaystyle 
	= -(\bar{\beta} - q)/q$}
\end{array}
\]
and, similarly,
\[
\liminf_{n \rightarrow \infty} \frac{1}{\log n} \log \mathbb{P}\left(\|S_{n}\| > s n^{1/p} \right) \leq  -(\underline{\beta} - q)/q.
\]
Letting $q \searrow p$, we have, for all $s > 0$,
 \begin{equation}
 \limsup_{n \to \infty} \frac{1}{\log n} \log \mathbb{P}\left(\left\|S_{n} \right\| > s n^{1/p} \right) 
 \leq - (\bar{\beta} - p)/p
 \end{equation}
 and 
 \begin{equation}
 \liminf_{n \to \infty} \frac{1}{\log n} \log \mathbb{P}\left(\left\|S_{n} \right\| > s n^{1/p} \right) 
 \leq - (\underline{\beta} - p)/p.
 \end{equation}
 Thus, (2.4) follows from (4.1) and (4.6) and, (2.5) follows from (4.2) and (4.7). The proof
 of Theorem 2.1 is complete. ~$\Box$
 
 \vskip 0.3cm

{\it Proof of Theorem 2.2}~~Since $\{X, X_{n}; n \geq 1\}$ is a sequence of i.i.d. real-valued random variables, we have for
$0 < p < 2$, (2.3) and (2.9) are equivalent and hence, Theorem 2.2 follows from Theorem 2.1 immediately.~ $\Box$

\vskip 0.3cm

{\it Proof of Theorem 2.3}~~For $0 < p < 1$, since $\{X, X_{n}; n \geq 1\}$ is a sequence of non-negative i.i.d. random variables, 
we have, for all $s > 0$,
\[
\mathbb{P}\left(S_{n} > sn^{1/p} \right) = \mathbb{P}\left(\left|S_{n}\right| > sn^{1/p} \right),~~n \geq 1,
\]
and hence, the conclusions of Theorem 2.3 (i) follow from (2.11) and Theorem 2.2. 

For $1 \leq p < 2$, since $X$ is non-negative random variable, under condition (2.12), we have
\begin{equation}
\frac{S_{n, \mu}}{n^{1/p}} \rightarrow_{\mathbb{P}} 0,
\end{equation}
where $X_{\mu} = X - \mu$, $X_{n, \mu} = X_{n}- \mu$, and $S_{n, \mu} = X_{1,\mu} + \cdots + X_{n,\mu}$, $n \geq 1$. 
It is easy to see that
\[
\limsup_{t \rightarrow \infty} \frac{1}{t} \log \mathbb{P}(\log|X - \mu| > t) 
=  \limsup_{t \rightarrow \infty} \frac{1}{t} \log \mathbb{P}(\log X > t) 
= - \bar{\alpha}
\]
and
\[
\liminf_{t \rightarrow \infty} \frac{1}{t} \log \mathbb{P}(\log|X - \mu| > t) 
=  \liminf_{t \rightarrow \infty} \frac{1}{t} \log \mathbb{P}(\log X > t) 
= - \underline{\alpha}.
\]
Thus, it follows from (4.8) and Theorem 2.2 that, for all $s > 0$,
\begin{equation}
\limsup_{n \rightarrow \infty} \frac{1}{\log n} 
\log \mathbb{P} \left(S_{n} > n \mu + sn^{1/p} \right) \leq \limsup_{n \rightarrow \infty} \frac{1}{\log n} 
\log \mathbb{P} \left(\left|S_{n,\mu} \right| > sn^{1/p} \right) = -\left(\bar{\alpha} - p \right)/p
\end{equation}
and
\begin{equation}
\liminf_{n \rightarrow \infty} \frac{1}{\log n} 
\log \mathbb{P} \left(S_{n} > n \mu + sn^{1/p} \right) \leq \liminf_{n \rightarrow \infty} \frac{1}{\log n} 
\log \mathbb{P} \left(\left|S_{n,\mu} \right| > sn^{1/p} \right) = -\left(\underline{\alpha} - p \right)/p.
\end{equation}
Write $S_{0, \mu} = 0$ and
\[
m_{n} = \max_{0 \leq k \leq n} m\left(S_{k, \mu} \right), ~n \geq 1.
\]
Since $X_{1, \mu}, ..., X_{n, \mu}$ are i.i.d. random variables, we have
\[
\max_{0 \leq k \leq n} m\left(S_{n, \mu} - S_{k, \mu} \right) = m_{n},~~n \geq 1.
\]
Hence, from (4.8), we have
\[
\lim_{n \rightarrow \infty} \frac{m_{n}}{n^{1/p}} = 0,
\] 
and furthermore, for given $s > 0$, by Lemma 3.7, for all sufficiently large $n$,
\begin{equation}
\begin{array}{lll}
\mbox{$\displaystyle 
	\mathbb{P} \left(\max_{1 \leq k \leq n} X_{k} > 3sn^{1/p}\right)$}
&\leq &
\mbox{$\displaystyle 
	\mathbb{P} \left(\max_{1 \leq k \leq n} X_{k, \mu} > 2sn^{1/p}\right)$}\\
&&\\
& \leq & 
\mbox{$\displaystyle 
	\mathbb{P} \left(\max_{1 \leq k \leq n} X_{k, \mu} > sn^{1/p} - 2m_{n} \right)$}\\
&&\\
& \leq & 
\mbox{$\displaystyle 
	\mathbb{P} \left(\max_{1 \leq k \leq n} \left(X_{k, \mu} 
	+ m\left(S_{k-1, \mu}\right) \right) > sn^{1/p} - m_{n} \right)$}\\	
&&\\
& \leq & 
\mbox{$\displaystyle 
	2 \mathbb{P} \left(\max_{1 \leq k \leq n} S_{k, \mu} > sn^{1/p} - m_{n} \right)$}\\	
&&\\
& \leq & 
\mbox{$\displaystyle 
	2 \mathbb{P} \left(\max_{1 \leq k \leq n} \left(S_{k, \mu} 
	+ m\left(S_{n, \mu} - S_{k, \mu}\right) \right) > sn^{1/p}\right)$}\\	
&&\\
& \leq & 
\mbox{$\displaystyle 
	4 \mathbb{P}\left(S_{n, \mu} > sn^{1/p} \right)$.}
\end{array}
\end{equation}
Since, under condition (2.12), $p \leq \bar{\alpha} \leq \underline{\alpha}$, by Lemma 3.4, (4.11) ensures that, for all $s > 0$,
\begin{equation}
\begin{array}{lll}
\mbox{$\displaystyle 
\limsup_{n \rightarrow \infty} \frac{1}{\log n} \log \mathbb{P} \left(S_{n} > n \mu + sn^{1/p} \right)$}
& = &
\mbox{$\displaystyle \limsup_{n \rightarrow \infty} \frac{1}{\log n} \log \mathbb{P}\left(S_{n, \mu} > sn^{1/p} \right)$}\\
&&\\
& \geq &
\mbox{$\displaystyle \limsup_{n \rightarrow \infty} 
	\frac{1}{\log n} \log \mathbb{P} \left(\max_{1 \leq k \leq n} X_{k} > 3sn^{1/p}\right)$}\\
&&\\
& = &
\mbox{$\displaystyle - \left(\bar{\alpha} - p \right)/p$}
\end{array}
\end{equation}
and, similarly,
\begin{equation}
\liminf_{n \rightarrow \infty} \frac{1}{\log n} 
\log \mathbb{P} \left(S_{n} > n \mu + sn^{1/p} \right) \geq  - \left(\underline{\alpha} - p \right)/p.
\end{equation}
Thus, (2.13) follows from (4.9) and (4.12) and, (2.14) follows from (4.10) and (4.13). This completes the proof
of Theorem 2.3. ~$\Box$

\vskip 0.5cm

\noindent {\bf Declarations}

\vskip 0.3cm

\noindent
{\bf Funding}\\
The research of Deli Li was partially supported by a grant from the Natural Sciences 
and Engineering Research Council of Canada (grant \#: RGPIN-2019-06065) and the research 
of Yu Miao was partially supported by a grant from the National Natural Science Foundation
of China (grant \#: NSFC-11971154).

\vskip 0.3cm

\noindent 
{\bf  Conflicts of interests/Competing interests}\\
The authors declare that they have no conflicts of interest.

\vskip 0.5cm

{\bf References}

\begin{enumerate}
	
	\item Azlarov, T. A., Volodin, N. A.: Laws of large numbers for identically distributed 
	Banach-space valued random variables. Teor. Veroyatnost. i Primenen. {\bf 26}, ~584-590 
	(1981), in Russian. English translation in Theory Probab. Appl. {\bf 26}, 573-580 (1981).
	
	\item Bahadur, R. R., Zabell, S. L.: Large deviations of the sample mean in general vector spaces.
	Ann. Probab. {\bf 7}, 587-621 (1979). 
	
	\item Chernoff, H.: A measure of asymptotic efficiency for tests of a hypothesis based on the sum
	of observations. Ann. Math. Statistics {\bf 23}, 493-507 (1952).
	
	\item Cram\'{e}r, H.: Sur un nouveau th\'{e}or\`{e}me-limite de la théorie des probabilit\'{e}s. 
	Actualit\'{e}s Sci. Indust. {\bf 736}, 5-23 (1938).
		
	\item de Acosta, A.: Inequalities for {\it B}-valued random vectors with applications 
	to the law of large numbers. Ann. Probab. {\bf 9}, 157-161 (1981).
	
	\item Donsker, M. D., Varadhan, S. R. S.: Asymptotic evaluation of certain Markov process
	expectations for large time. III. Comm. Pure Appl. Math. {\bf 29}, 389-461 (1976).
	
	\item Feller, W.: An Introduction to Probability Theory and its Applications, Vol. 2, 2nd ed., 
	Wiley, New York (1971).
	
	\item Gantert, N.: A note on logarithmic tail asymptotics and mixing. Statist. Probab. Lett. 
	{\bf 49}, 113-118 (2000).
	
	\item Gnedenko, B. V., Kolmogorov, A. N.: Limit Distributions for Sums of Independent Random
	Variables, 2nd ed., Addison-Wesley, Reading, MA (1968).
	
	\item Hoffmann-J{\o}rgensen, J.: Sums of independent Banach space valued random variables. 
	Studia Math. {\bf 52}, 159-186 (1974).
	
	\item Hu, Y. J., Nyrhinen, H.: Large deviations view points for heavy-tailed random walks. 
	J. Theoret. Probab. {\bf 17}, 761-768 (2004).
	
	\item Kolmogoroff, A.: Sur la loi forte des grands nombres. C. R. Acad. Sci. Paris Ser. Math. 
	{\bf 191}, 910-912 (1930).
	
	\item Ledoux, M., Talagrand, M.: Probability in Banach Spaces:	Isoperimetry and Processes. 
	Springer-Verlag, Berlin (1991).
	
	\item Li, D., Liang, H.-Y., Rosalsky, A.: A probability inequality for sums of independent
	Banach space valued random variables. Stochastics {\bf 90}, 214-223 (2018). 
	
	\item Li, D., Rao, M. B., Jiang, T., Wang, X.: Complete convergence and almost sure 
	convergence of weighted sums of random variables. J. Theoret. Probab. {\bf 8}, 49-76 (1995).
	
	\item Li, D., Rosalsky, A.: Precise lim sup behavior of probabilities of large deviations 
	for sums of i.i.d. random variables. Internat. J. Math. \& Math. Sci. {\bf 2004}, 3565-3576 (2004).
	
	\item Marcinkiewicz, J., Zygmund, A.: Sur les fonctions ind\'{e}pendantes. Fund. Math. {\bf 29}, 60-90 (1937).
	
	\item Marcus, M. B., Woyczy\'{n}ski, W. A.: Stable measures and	central limit theorems 
	in spaces of stable type. Trans. Amer. Math. Soc. {\bf 251}, 71-102 (1979).
		
	\item Miao, Y., Xue, T. Y., Wang, K., Zhao, F. F.: Large deviations for dependent heavy tailed random variables. 
	J. Korean Statist. Soc. {\bf 41}, 235-245 (2012).
	
	\item Mourier, E.: El\'{e}ments al\'{e}atoires dans un espace de Banach. Ann. Inst. H. Poincar\'{e} 
	{\bf 13}, 161-244 (1953).
	
	\item Petrov, V. V.: A generalization of an inequality of L\'{e}vy. Teor Veroyatnost. i Primenen. 
	{\bf 20}, 140-144 (1975) (Russian), translated in Theory Probab. Appl. {\bf 20}, 141-145 (1975).
	
	\item Petrov, V. V.: Limit Theorems of Probability Theory: Sequences of Independent Random Variables. 
	Clarendon Press, Oxford (1995).
	
	\item Wang, X. J.: Upper and lower bounds of large deviations for some dependent sequences. Acta Math.
	Hungar.  {\bf 153}, 490-508 (2017).
	
\end{enumerate}

\end{document}